\documentclass{llncs}

\usepackage[utf8]{inputenc}
\usepackage[T1]{fontenc}
\usepackage[brazilian]{babel}
\usepackage{csquotes}
\usepackage{hyphenat}
\usepackage{amssymb}
\usepackage{amsmath}
\usepackage{array}

\usepackage[locale=FR]{siunitx}

\usepackage{graphicx}
\graphicspath{ {images/} }

\usepackage{listings}

\begin{document}

\title{Nontrivial Solutions of Dirac-Laplace Equation on Compact Spin Manifolds}
%

\author{Xu Yang* and Lei Xian}
%
%
%
\institute{School of Mathematics , Yunnan Normal University\\
\email{yangxu02465129@163.com}}

\maketitle              
\renewcommand\refname{References}
\renewcommand{\abstractname}{Abstract}
\begin{abstract}
We apply the Fountain theorem
to a class of nonlinear Dirac-Laplace equation on  compact spin manifold.  We show the
standard Am brosetti-Rabinowitz condition can be replaced by a more natural super-quadratic condition that is enough to obtain the Cerami condition under certain conditions.
Multiple solutions of nonlinear Dirac-Laplace equation are obtained in this note.
\end{abstract}
\section{Introduction and main results}
Motivated by quantum physics, Esteban-S\'er\'e \cite{EsSe1,EsSe2} studied
 existence and multiplicity of  solutions of nonlinear Dirac
equations on $\mathbb{R}^3$. In this case, a large number of existence of solution has been obtained.
Dirac operators on compact spin manifolds play prominent role in the geometry and mathematical physics, such as the generalized Weierstrass representation of the surface in three manifolds \cite{Fri} and the supersymmetric nonlinear sigma model in quantum field theory \cite{CJLW1,CJLW2}.
 For the nonlinear Dirac equations on a general compact spin manifold, some results were recently obtained,
see Amann \cite{Amm}, Isobe \cite{Iso1,Iso2},  Maalaoui and V. Martino \cite{Maa,MaM}
, Gong and Lu \cite{GoL1,GoL2} and \cite{YJL,Yan} by authors. In addition, be different with these existing works, Ding and Li \cite{DiL} studied  a class of boundary value problem on a compact spin manifold $M$ with smooth boundary. The problem is a general relativistic model of confined particles by means of nonlinear Dirac fields on $M$. In this paper, we are concerned with a nonlinear square Dirac equation on compact spin manifolds, and deal with some existence results.

Assume $m\ge 2$ is an integer.
 Given an $m$-dimensional compact oriented Riemannian manifold
  $(M,g)$ equipped with a spin structure $\rho:P_{Spin(M)}\rightarrow P_{SO(M)}$,  let $\Sigma M=\Sigma(M,g)=P_{Spin(M)} \times_{\sigma}\Sigma_m$ denote
the complex spinor bundle on $M$, which is a complex vector bundle of rank $2^{[m/2]}$ endowed with the spinorial Levi-Civita connection $\nabla$ and a pointwise Hermitian scalar product $\langle\cdot,\cdot\rangle$.  Write the point of $\Sigma M$ as
$(x, \psi)$, where $x\in M$ and $\psi\in \Sigma_x M$. Laplace operator on spinors $\Delta:C^\infty(M,\Sigma M)\rightarrow C^\infty(M,\Sigma M)$ is defined by
$$
\Delta\psi=-\sum_{i=1}^n(\nabla_{e_i}\nabla_{e_i}-\nabla_{\nabla_{e_i}e_i})\psi,
$$
 it is independent of the choice of local orthonormal frame  and  is a second order differential elliptic operator.

 The Dirac operator $D$ is an elliptic differential operator of order one,
 $$
 D=D_g:C^{\infty}(M,\Sigma M)\rightarrow C^{\infty}(M,\Sigma M),
 $$
 locally given by $D\psi=\sum^m_{i=1}e_i\cdot\nabla_{e_i}\psi$ for $\psi\in
C^{\infty}(M,\Sigma M)$ and a local $g$-orthonormal frame
$\{e_i\}^m_{i=1}$ of the tangent bundle $TM$. It is an unbounded essential self-adjoint operator in $L^2(M,\Sigma M)$ with the domain $C^\infty(M,\Sigma M)$,
and its spectrum consists of an unbounded sequence of real numbers.
these mean that the closure of $D$, is a self-adjoint operator in $L^2(M,\Sigma M)$  with the domain $H^1(M,\Sigma M)$.

The square $D^2$ of the Dirac operator as well as that of the Laplace operator $\Delta$ is second order differential operator, which is called Dirac-Laplace operator. On the Dirac operator and square Dirac operator under boundary condition has been discussed a lot in the literature,c.f.\cite{DiL,FaSc,Fri,WaK,WoC}. the connection between $D^2$ and $\Delta$ is given by Schr\"{o}dinger-Lichnerowicz formula
\begin{equation*}
  D^2=\Delta+\frac{R}{4}\text{Id},
\end{equation*}
where $R$ is the scalar curvature of $(M,g)$.

For a fiber preserving nonlinear map $h:\Sigma M\rightarrow \Sigma M$,
we consider the following nonlinear Dirac-Laplace equations
\begin{eqnarray}\label{e:1.1}
D^2\psi(x)= h(x,\psi(x)) \quad\hbox{on}\;\ M,
\end{eqnarray}
where $\psi \in C^\infty(M,\Sigma M)$ is a spinor. We assume $h$ has a $\psi$-potential, namely there exist real valued function $H:\Sigma M\rightarrow \mathbb{R}$ such that $H_\psi=h$, then equation (\ref{e:1.1}) has a variational structure.
 It is  the Euler-Lagrange equation of the functional
\begin{eqnarray}\label{e:1.2}
\mathfrak{L}(\psi)=\frac{1}{2}\int_M\langle D\psi,D\psi\rangle dx
-\int_M H(x,\psi(x))dx
\end{eqnarray}
on the Sobolev space $H^1(M,\Sigma M)$, see\cite{Ada} and \cite{Iso1}. Here  $dx$ is the Riemannian volume measure on $M$ with respect to the metric $g$, and $\langle\cdot,\cdot\rangle$ is the compatible metric on $\Sigma M$.

The functional(\ref{e:1.2}) is strongly indefinite in the sense that  the Dirac-Laplace operator $D^2$ takes  infinitely
 many positive eigenvalues. This kind of problem similar to elliptic equation $\Delta \psi= h(x,\psi)$ on $\mathbb{R}^n$. In\cite{AmR}, Ambrosetti and Rabinowitz established the Mountain Pass Theorem which successfully  applied to the elliptic equation on $\mathbb{R}^n$. It plays a key role not only in establishing the mountain pass geometry of the functional but also in obtaining Palais-Smale condition. It has been an open question that whether the more natural superlinear  condition is enough to obtain the A-R condition for variational problem. In\cite{LiW}, the existence of solutions is given on $\mathbb{R}^n$  without Ambrosetti-Rabinowitz condition. Recently, Isobe\cite{Iso1} consider Dirac equations on compact spin manifold via Galerkin type approximations and linking arguments, the existence of infinitely many solutions is obtained under  Ambrosetti-Rabinowitz condition.
 As showed by  Liu and Li\cite{LiL}, the Fountain theorem under Cerami condition also provided a powerful method to deal with strongly  indefinite variational problems.
In this note our aim is to study (\ref{e:1.1}) with this method.
 Our techniques work for the existence of solutions to Dirac-Laplace equation under Cerami condition replace Palais-Smale condition.

For the nonlinearity $H$, we make the following hypotheses:

$H$ is $\alpha$-H\"{o}lder continuous in the direction of the base, where $0<\alpha<1$.

\noindent{(${\bf H}_1$)} There exist a constant $C>0$ such that
\begin{eqnarray}
  \left| H_{\psi}(x,\psi)\right|&\leq& C\left(1+|\psi|^{q-1}\right)\label{e:1.3}
\end{eqnarray}
for any $(x,\psi)\in \Sigma M$, where $2<q< 2^{\ast}$, $2^{\ast}=\frac{2m}{m-2}$ if $m\geq 3$ and $2^{\ast}=\infty$ if $m=2$.
\\
\noindent{(${\bf H}_2$)} $\frac{\langle H_\psi(x,\psi),\psi\rangle}{\mid\psi\mid^2}\rightarrow +\infty$ as $|\psi|\rightarrow +\infty$ uniformly for $x\in M$.
\\
\noindent{(${\bf H}_3$)} There exist $\theta\geq 1$ such that  $\theta \hbar(x,\psi)\geq \hbar(x,\mu \psi)$, where $(x,\psi)\in \Sigma M$, $\mu \in [0,1]$ and $ \hbar(x,\psi)=\langle H(x,\psi),\psi\rangle-2H(x,\psi)$.
\\
\noindent{(${\bf H}_4$)} $H(x,-\psi)=H(x,\psi)$ for any $(x,\psi)\in \Sigma M$.

Our main result is as follow.
\begin{theorem}\label{the:1.1}
If the above $H$ satisfies $({\bf H}_1)-({\bf H}_4)$,
then there exist a sequence of solutions $\{\psi_k\}_{k=1}^\infty\subset C^2(M,\Sigma M)$ to (\ref{e:1.1}) with $\mathfrak{L}(\psi_k)\rightarrow \infty$ as $k\rightarrow \infty$.
\end{theorem}

{\bf Notice:} The standard Ambrosetti-Rabinowitz superlinear condition has appeared In Isobe\cite{Iso1} as following:

(AR) There exist $\mu>2$ and $R_1>0$ such that
$$
0<\mu H(x,\psi)\leq \langle H_\psi(x,\psi),\psi\rangle
$$
for any $(x,\psi)\in \Sigma M$ with $\mid \psi\mid\geq R_1$. It implies there is $C>0$ such that $F(x,\psi)\geq C\mid \psi\mid^\mu$ for $\mid \psi\mid$ large. Isobe\cite{Iso1} obtained existence and multiplicity of solutions to nonlinear Dirac equations on compact spin manifolds. However,
in this paper, we use a more natural and weaker superlinear condition is that $(\bf{H_2})$:
$$
\frac{\langle H_\psi(x,\psi),\psi\rangle}{\mid\psi\mid^2}\rightarrow +\infty
$$
as $|\psi|\rightarrow +\infty$ uniformly for $x\in M$.

In this paper, condition $(\bf{H_3})$ refers to the condition of Liu and Li\cite{LiL}.  Notice that $H(x,\psi)=H(x)|\psi|^p$ satisfies these conditions$(\bf{H_1})-(\bf{H_4})$. We also remark that the H\"{o}ldre condition of $H$ is only used to prove the $C^2$-regularity of $\psi$.

Finally, as a concluding remark we point out that
in $\mathbb{R}^n$, the scalar curvature vanished,i.e. $R=0$. Then $D^2=\triangle $ with $\Delta=-\sum\frac{\partial^2}{\partial x_i^2}$. The Dirac-Laplace equations (\ref{e:1.1})
reduces to the following the elliptic equation
\begin{equation*}
  \Delta \psi= h(x,\psi).
\end{equation*}
For this elliptic equation, the existence of solutions has been obtained under assumption Ambrosetti-Rabinowitz condition and different additional conditions.

\section{Preliminaries}
According to \cite{Fri,LaM}, we next give the definition of Laplace operator on spinors:
\begin{definition}\label{def:2.1}
If $\psi\in C^\infty(M,\Sigma M)$, then $\Delta:C^\infty(M,\Sigma M)\rightarrow C^\infty(M,\Sigma M)$ is defined by
$$
\Delta\psi=-\sum_{i=1}^n(\nabla_{e_i}\nabla_{e_i}-\nabla_{\nabla_{e_i}e_i})\psi.
$$
\end{definition}
Using Stokes'theorem we have
$$
\int_M\langle \Delta\psi,\psi\rangle=\int_M\langle \nabla\psi,\nabla\psi\rangle=\int_M\langle\psi,\Delta\psi\rangle
$$
for any $\psi\in C^\infty(M,\Sigma M)$. Here, $\langle \nabla\psi,\nabla\psi\rangle$ is the scalar product on $1$-forms,i.e.
$\langle \nabla\psi,\nabla\psi\rangle=\sum\limits_{i=1}^n\langle \nabla_{e_i}\psi,\nabla_{e_i}\psi\rangle$, where $\{e_i\}_{i=1}^\infty$ is a local
orthonormal frame on the spin manifold.

The first Sobolev norm \cite{Fri} of a smooth spinor field $\psi\in C^\infty(M,\Sigma M)$ is given by
$$
\parallel \psi\parallel_{H^1}^2=\parallel\psi\parallel_2^2+\parallel\nabla\psi\parallel_2^2
$$
and the corresponding Sobolev space $H^1(M,\Sigma M)$ is the completion of $\psi\in C^\infty(M,\Sigma M)$ with respect to this norm, where
$\parallel\cdot\parallel_2^2:=\int_M\mid\psi\mid^2dx $.

The Dirac operator D acts on spinors on $M$, $D:C^\infty(M,\Sigma M)\rightarrow C^\infty(M,\Sigma M)$ in the following definition:
\begin{definition}\label{def:2.2}
The composition
\begin{eqnarray*}
  D &=& c\circ\bigtriangledown:C^\infty(M,\Sigma M)\rightarrow C^\infty(M,T^*M\oplus\Sigma M) \\
   &=& C^\infty(M,TM\oplus\Sigma M)\rightarrow C^\infty(M,\Sigma M)
\end{eqnarray*}
is called the Dirac operator.
\end{definition}
It is locally given by $D\psi=\sum^m_{i=1}e_i\cdot\nabla_{e_i}\psi$ for $\psi\in
C^{\infty}(M,\Sigma M)$. It is a first order elliptic operator in $L^2(M,\Sigma M)$ with the domain $C^\infty(M,\Sigma M)$. For a vector $X\in TM$, its symbol $\sigma(D)(X):\Sigma M\rightarrow\Sigma M$ is given by Clifford multiplication:$\sigma(D)(X)(\psi)=X\cdot\psi$. By the Dirac operator definition, we know the Dirac operator is a symmetric operator on compact spin manifold,i.e.
$$
\int_M\langle D\psi,\phi\rangle dx=\int_M\langle \psi,D\phi\rangle dx,
$$
where the  smooth spinor field $\psi,\phi\in C^\infty(M,\Sigma M)$.

Another Sobolev norm is given  by
$$
\parallel \psi\parallel_{1,2}^2=\parallel\psi\parallel_2^2+\parallel D\psi\parallel_2^2.
$$
Thus we will obtain the next Lemma:
\begin{lemma}\label{lem:2.3}
The norm $\parallel \psi\parallel_{H^1}$ and the norm $\parallel \psi\parallel_{1,2}$ are equivalent norms.
\end{lemma}
proof: Firstly, from the estimate below for the norm of $D\psi$:
\begin{eqnarray*}
\parallel D\psi\parallel_2^2 &=& \sum\limits_{i,j=1}^n\int_M\langle e_i\cdot\nabla_{e_i}\psi,e_j\cdot\nabla_{e_j}\psi\rangle dx\leq
\sum\limits_{i,j=1}^n\int_M\mid\nabla_{e_i}\psi\mid\mid\nabla_{e_j}\psi\mid dx \\
   &\leq& \frac{1}{2}\sum\limits_{i,j=1}^n\int_M(\mid\nabla_{e_i}\psi\mid^2+\mid\nabla_{e_j}\psi\mid^2)dx=n\parallel \nabla\psi\parallel_2^2\leq n\parallel\psi\parallel_{H^1}^2.
\end{eqnarray*}
This shows the Dirac operator $D:H^1(M,\Sigma M)\rightarrow L^2(M,\Sigma M)$ is a continuous operator and
\begin{equation}\label{e:2.1}
  \parallel \nabla\psi\parallel_2^2\geq\frac{1}{n}\parallel D\psi\parallel_2^2.
\end{equation}
On the other hand, by Schr\"{o}dingger-Lichnerowicz formula, rewrite $\parallel D\psi\parallel_2^2$, we have as follows:
\begin{equation}\label{e:2.2}
 \parallel D\psi\parallel_2^2=\int_M\langle D^2\psi,\psi\rangle dx=\parallel \nabla\psi\parallel_2^2+\frac{R}{4}\parallel\psi\parallel_2^2.
\end{equation}
Set $R_{min}=\min\{R(x):x\in M\}$ is the minimum of the scalar curvature and $R_{max}$ its maximum. By the (\ref{e:2.2}), we have
\begin{equation}\label{e:2.3}
  \parallel\psi\parallel_{H^1}^2+(\frac{R_{min}}{4}-1)\parallel\psi\parallel_2^2\leq \parallel D\psi\parallel_2^2\leq\parallel\psi\parallel_{H^1}^2+(\frac{R_{max}}{4}-1)\parallel\psi\parallel_2^2.
\end{equation}
From (\ref{e:2.1}) and (\ref{e:2.3}), we obtain
$$
\frac{1}{n}(\parallel\psi\parallel_2^2+\parallel D\psi\parallel_2^2)\leq\parallel\psi\parallel_{H^1}^2\leq\parallel D\psi\parallel_2^2+(1-\frac{R_{min}}{4})\parallel\psi\parallel_2^2.
$$
Thus $\parallel \psi\parallel_{H^1}$ and $\parallel \psi\parallel_{1,2}$ are equivalent norms. In other words,
the Sobolev space $H^1(M,\Sigma M)$ can be defined as the completion of the $C^\infty(M,\Sigma M)$ with respect to the norm $\parallel \psi\parallel_{1,2}$. In the following, we will use the norm $\parallel \psi\parallel_{1,2}$ respect to the $H^1(M,\Sigma M)$.

The operator $D$ and $D^2$ are essentially self-adjoint in $L^2(M,\Sigma M)$. The eigenspinor $\psi_k$
of $D$(with eigenvalue $\lambda_k$) are also eigenspinor of $D^2$ (with eigenvalue $\lambda_k^2$),i.e.
$$
spec(D^2)=\{\lambda_k^2| \lambda_k\in spec(D)\},
$$
where $spec(D)$ is expressed as a set of all eigenvalues of $D$ on the $H^1(M,\Sigma M)$.
\par
In addition, the following Lemma c.f.\cite{Fri} is showed that the kernels of $D$ and $D^2$ in $L^2(M,\Sigma M)$ is coincide,i.e. $ker(D)=ker(D^2)$.
\begin{lemma}\label{lem:2.4}
Let $M$ be a compact spin manifold and $\psi \in C^\infty(M,\Sigma M)$. Then for any number $t>0$,
$$
\parallel D\psi\parallel_2^2\leq t\parallel D^2\psi\parallel_2^2+\frac{1}{t}\parallel\psi\parallel_2^2.
$$
\end{lemma}
In the next, we will prove the conclusion $ker(D)=ker(D^2)$.

Let $\psi\in L^2(M,\Sigma M)$ satisfy $D^2\psi=0$. By the regularity theorem for solution of elliptic differential equations  we first conclude that $\psi$ is smooth. By the Lemma(\ref{lem:2.4}), we obtain
$$
\parallel D\psi\parallel_2^2\leq t\parallel D^2\psi\parallel_2^2+\frac{1}{t}\parallel\psi\parallel_2^2=\frac{1}{t}\parallel\psi\parallel_2^2.
$$
Since $\parallel\psi\parallel_2^2<\infty$ imply for $t\rightarrow \infty$ that $D\psi\equiv 0$. Hence, $ker(D)=ker(D^2)$.

We introduce notation which will be frequently used throughout this paper. The complete orthonormal basis $\{\psi_k\}$ of $L^2(M,\Sigma M)$ consisting  of the eigenspinors of $D^2$ is decomposed into two parts: $\{\psi_k\}=\{\psi_k^+\}_{k=1}^\infty\cup \{\psi_k^0\}_{k=1}^\kappa$, where
where $D^2\psi_k^+=\lambda_k^2\psi_k^+$ with $\lambda_k^2>0$; $D^2\psi_k^0=0$ and $\kappa=dim (ker D)=dim (ker D^2)<\infty$. By the elliptic regularity, we have any eigenspinor $\psi_k\in C^\infty(M,\Sigma M)$. Hence, we also set $H^+$ be the subspace spanned by eigenspinors $\{\psi_k^+\}_{k=1}^\infty$ with positive eigenvalues $\{\lambda_k^2\}_{k=1}^\infty$, and $H^0$ the nullspace of $D^2$. Then we have the decomposition of the Hilbert space $H^1(M,\Sigma M)$:
$$
H^1(M,\Sigma M)=H^0\oplus H^+.
$$
where $H^0:= \overline{span\{\psi_k^0\}_{k=1}^\kappa}$, $H^+:= \overline{span\{\psi_k^+\}_{k=1}^\infty}$.

Under the condition $(H_1)$, by  Sobolev embedding theorem, it is easily checked that $\mathfrak{L}$ is a $C^1$ functional on $H^1(M,\Sigma M)$ and we obtain the following proposition:

\begin{proposition}\label{prop:2.5}
Under the condition (${\bf H}_1$) the functional  $\mathcal{H}:H^{1}(M,\Sigma M)\to \mathbb{R}$
 defined by
 \begin{equation} \label{e:2.4}
\mathcal{H}(\psi)=\int_MH(x,\psi(x))dx,
\end{equation}
 is of class $C^1$, and  at each $\psi\in H^{1}(M,\Sigma M)$ derivations $\mathcal{H}^\prime(\psi)$
 was given by
\begin{equation}\label{e:2.5}
\mathcal{H}^\prime(\psi)\xi=\int_M\langle H_\psi(x,\psi(x)),\xi(x)\rangle dx\quad\forall \xi\in H^{1}(M,\Sigma M).
\end{equation}
\end{proposition}
In the view of the calculus of variations, the weak solutions to the problem (\ref{e:1.1}) are obtained as critical points of the following Euler-Lagrange functional
\begin{eqnarray*}
\mathfrak{L}(\psi)=\frac{1}{2}\int_M\langle D\psi,D\psi\rangle dx
-\int_M H(x,\psi(x))dx
\end{eqnarray*}
By the standard elliptic regularity theory, such a weak solution is in fact $C^2$ and a classical solution to problem (\ref{e:1.1}).

\section{The Cerami condition for $\mathfrak{L}$}
Let $F$ be a $C^1$ functional  on a Banach space $E$, $c\in \mathbb{R}$. Recall that
 a sequence $\{u_{n}\}\subset E$ is called a Cerami condition if
 \\
 \noindent{(i)} There exists a convergent subsequence  for  any boundedness sequence $\{u_n\}$ with $F(u_{n})\rightarrow c$  and
 $\parallel dF(u_{n})\parallel_{E^{*}}\rightarrow 0$ as $n\rightarrow \infty$.
 \\
 \noindent{(ii)}There exists $\delta,r,\beta>0$ such that
 $$
 \parallel dF(u_n)\parallel_{E^{*}}\parallel u_n\parallel\geq \beta
 $$
 for any $u_n\in F^{-1}[c-\delta,c+\delta]$ with $\parallel u_n\parallel\geq r$.
 \\
In this section we prove the  Cerami condition for $\mathfrak{L}$.

\begin{lemma}\label{lem:3.1}
Suppose $H$ satisfies $({\bf H}_1)$,$({\bf H}_2)$ and $({\bf H}_3)$. Then for any $c\in \mathbb{R}$, $\mathfrak{L}$ satisfies the Cerami condition with respect to $H^1(M,\Sigma M)$.
\end{lemma}
\noindent{\bf Proof}.\quad    We prove first Cerami conditon (i) is satisfied. Let $\psi_n$ be bounded in $H^1(M,\Sigma M)$. By the Sobolev embedding theorem, we have the compact embedding
$$
H^1(M,\Sigma M)\hookrightarrow L^p(M,\Sigma M)
$$
for $1\leq p<2^*$. Passing to a subsequence, we may assume that for some $\psi\in H^1(M,\Sigma M)$,$\psi_n\rightharpoonup\psi$ weakly in  $H^1(M,\Sigma M)$
and $\psi_n\rightarrow \psi$ strongly in $L^p(M,\Sigma M)$ for any $1\leq p<2^*$.
Setting $\psi_n=\psi_n^0+\psi_n^+$ and $\psi=\psi^0+\psi^+$ according to the $H^1(M,\Sigma M)=H^0\bigoplus H^+$.

since $\text{dim} (H^0)<\infty$, $\parallel\parallel_{1,2}$ and $\parallel \parallel_2$ on $H^0$ are equivalent. Hence, we have
\begin{equation}\label{e:3.1}
  \parallel\psi_n^0-\psi^0\parallel_{1,2}\leq C\parallel\psi_n^0-\psi^0\parallel_2\rightarrow 0
\end{equation}
as $n\rightarrow \infty$.
\\
From the condition $(H_1)$, we get
\begin{eqnarray}\label{e:3.2}
  &&\parallel\psi_n^+-\psi^+\parallel^2_{1,2} = \parallel D\psi_n^+-D\psi^+\parallel_2^2+\parallel\psi_n^+-\psi^+\parallel_2^2\nonumber \\
&\leq& \parallel D\psi_n^+-D\psi^+\parallel_2^2+(\lambda_1^+)^{-2}\parallel D\psi_n^+-D\psi^+\parallel_2^2\leq C\parallel D(\psi_n^+-\psi^+)\parallel_2^2\nonumber\\
&\leq& C\langle d\mathfrak{L}(\psi_n^+-\psi^+ ),\psi_n^+-\psi^+ \rangle+C\int_{M} \langle H_\psi(x,\psi_n^+-\psi^+),\psi_n^+-\psi^+ \rangle dx\nonumber \\
&\leq& C\parallel d\mathfrak{L}(\psi_n^+-\psi^+ )\parallel_{H^{1*}}\parallel \psi_n^+-\psi^+ \parallel_{1,2}+C\int_{M}  (1+\mid\psi_n^+-\psi^+\mid^{p-1})\mid\psi_n^+-\psi^+ \mid dx\nonumber \\
&\leq& o(1)+C(\parallel \psi_n^+-\psi^+\parallel_1+\parallel \psi_n^+-\psi^+\parallel_p^p)\rightarrow 0,
\end{eqnarray}
as $n\rightarrow\infty$.

(\ref{e:3.1}) and (\ref{e:3.2})  imply the subsequence of $\psi_n$ that converges
$\psi$ in $H^1(M,\Sigma M)$. So the Cerami condition (i) is verified.
\\
To prove the Cerami condition (ii), we assume that
\begin{equation}\label{e:3.3}
 \mathfrak{L}(\psi_n)\rightarrow c,\quad \parallel\psi_n\parallel_{1,2}\rightarrow \infty,\quad \parallel d\mathfrak{L}(\psi_n)\parallel_{H^{1*}}\parallel \psi_n\parallel_{1,2}\rightarrow 0
\end{equation}
as $n\rightarrow\infty$.
\\
From (\ref{e:3.3}) , we have
\begin{equation}\label{e:3.4}
  2\mathfrak{L}(\psi_n)-\langle d\mathfrak{L}(\psi_n),\psi_n\rangle
=\int_{M} \langle H_{\psi}(x,\psi_n),\psi_n\rangle dx-2\int_MH(x,\psi_n)dx\rightarrow 2c
\end{equation}
as $n\rightarrow \infty$.
\\
Set $\omega_n=\frac{\psi_n}{\parallel D\psi_n \parallel_2}$, which imply $\parallel \omega_n\parallel_{1,2}\leq C$.  By the $H^1(M,\Sigma M)$ is a reflexive space and Sobolev embedding theorem, there exist $\omega$ such that
\begin{eqnarray}\label{e:3.5}
  &&\omega_n\rightharpoonup \omega \quad\text{in } H^1(M,\Sigma M)\nonumber \\
  &&\omega_n\rightarrow \omega \quad\text{in } L^p(M,\Sigma M)\nonumber \\
  &&\omega_n(x)\rightarrow \omega(x)\quad a.e.\ x\in M
\end{eqnarray}
where $1\leq p<2^*$.

On the one hand, if $\parallel\omega(x)\parallel_2\equiv 0 $, then $\omega(x)=0$ a.e. on $M$.
\\
Define $I(t):[0,1]\rightarrow R$ such that
$$
I(t)=\mathfrak{L}(t\psi)=t^2\int_M |D\psi|^2dx-\int_MH(x,t\psi)dx.
$$
According to the techniques of Jeanjean\cite{Jea}, we define the sequence $\{t_n\}$ with $t_n\in [0,1]$ satisfying
\begin{equation}\label{e:3.6}
  I(t_n\psi_n)=\max\limits_{t\in [0,1]}\mathfrak{L}(t\psi_n).
\end{equation}
If, for $n\in \mathbb{N}$, $t_n$ is defined by $(\ref{e:3.6})$ is not unique, we choose any one of those values.
\\
For any $m>0$, set $\overline{\omega}_n(x)=\sqrt{2m}\omega_n(x)$, by $\omega_n(x)\rightarrow \omega(x)=0, a.e. x\in M$ and $(\bf{H}_1)$, it is easy to see that
\begin{equation}\label{e:3.7}
  \lim\limits_{n\rightarrow \infty}\int_MH(x,\overline{\omega}_n(x))dx=\int_M \lim\limits_{n\rightarrow \infty}H(x,\overline{\omega}_n(x))dx=0.
\end{equation}
Therefore, for enough $n$, according to $\parallel\psi_n\parallel_{1,2}\rightarrow \infty$, (\ref{e:3.6}) and (\ref{e:3.7}), we have
$$
\mathfrak{L}(t_n\psi_n)\geq\mathfrak{L}(\overline{\omega}_n(x))=2m-\int_MH(x,\overline{\omega}_n(x))dx\geq m.
$$
The above argument shows that $\mathfrak{L}(t_n\psi_n)\rightarrow \infty$ as $n\rightarrow \infty$.
\\
Since $\mathfrak{L}(0)=0$ and $\mathfrak{L}(\psi_n)\rightarrow c$, we have $0<t_n<1$. Hence for enough $n$, we conclude that
$$
d\mathfrak{L}(t_n\psi_n),t_n\psi_n\rangle=t_n\frac{d}{dt}\mid_{t=t_n}\mathfrak{L}(t\psi_n)=0.
$$
It implies that
\begin{equation}\label{e:3.8}
  \int_M|D(t_n\psi_n)|^2dx=\int_M\langle H_\psi(x,t_n\psi_n),t_n\psi_n\rangle dx.
\end{equation}
According to $0<t_n<1$ and $(\bf{H_3})$, there exist $\theta\geq1$ such that $\theta\hbar(x,\psi_n)\geq \hbar(x,t_n\psi_n)$. plugging (\ref{e:3.8}) into the following inequality, we obtain
\begin{eqnarray*}
  2\mathfrak{L}(\psi_n)-\langle  d\mathfrak{L}(\psi_n),\psi_n\rangle &=& \int_M \hbar(x,\psi_n)\geq \frac{1}{\theta}\hbar(x,t_n\psi_n)dx \\
   &=& \frac{1}{\theta}(\int_M \langle H_\psi(x,t_n\psi_n),t_n\psi_n\rangle dx-2\int_M H(x,t_n\psi_n))\\
   &=&\frac{2}{\theta}\mathfrak{L}(t_n\psi_n)\rightarrow \infty
\end{eqnarray*}
as $n\rightarrow \infty$. This is a contraction to (\ref{e:3.4}).
\\
On the other hand, if $\parallel\omega(x)\parallel_2\neq 0 $, it follows from$(\bf{H_2})$ and Fatou's Lemma that
$$
\int_M\frac{\langle H_\psi(x,\psi_n),\psi_n\rangle}{\mid \psi_n \mid^2}\mid \omega_n\mid^2dx\rightarrow \infty
$$
as $n\rightarrow \infty$.
But from(\ref{e:3.3}), we get
$$
\int_M\mid D\psi_n\mid^2 dx-\int_M \langle h(x,\psi_n),\psi_n\rangle dx=\langle d\mathfrak{L}(\psi_n),\psi_n\rangle=o(1)
$$
as $n\rightarrow \infty$. This is
$$
1-o(1)=\int_M\frac{\langle h(x,\psi_n),\psi_n\rangle}{\parallel D\psi_n \parallel_2^2}dx=\int_M\frac{\langle h(x,\psi_n),\psi_n\rangle}{\mid \psi_n \mid^2}\mid \omega_n\mid^2dx.
$$
Thus we get a contradiction. The proof is complete, $\mathfrak{L}$ satisfies the Cerami condition with respect to $H^1(M,\Sigma M)$.

\section{Proofs of Theorem~\ref{the:1.1}}
To obtain the theorem \ref{the:1.1}, we recall Fountain theorem for semi-definite functionals, see \cite{Wil} for the detailed exposition.
\\
Let  $X$ be a Banach space  with basis $\{ e_j\}_{j=1}^{\infty}$,i.e. $X=\overline{span\{ e_j\}_{j=1}^{\infty}}$. We set
\begin{equation*}
  Y_k:= \bigoplus\limits_{j=1}^{k}\mathbb{R}e_j,\quad Z_k:= \overline{\bigoplus\limits_{j=k}^{\infty}\mathbb{R}e_j},
\end{equation*}
We then have $X=Y_k\bigoplus Z_k$.

\begin{theorem}\label{the:4.1}
(Fountain theorem) Let $J \in C^1(X,\mathbb{R})$ is an even functional,i.e. $J(-u)=J(u)$ for all $u\in X$. $J$ satisfies  the Cerami condition.
 If for every $k\in \mathbb{N}$, there
exists $\rho_k>r_k>0$ such that
\begin{flalign*}
&(A_1): a_k:=\inf_{\substack{u\in Z_k\\ \parallel u\parallel= r_k}}J(u)\rightarrow \infty, k\rightarrow \infty; &\\
&(A_2): b_k:=\sup_{\substack{u\in Y_k\\ \parallel u\parallel=\rho_k}}J(u)\leq 0\
\end{flalign*}
then $J$ has an unbounded sequence of critical values.
\end{theorem}

We shall apply the fountain theorem for the functional $\mathfrak{L}$ on $H^1(M,\Sigma M)$. First of all,we prove:
\begin{lemma}\label{lem:4.2}
Define $\beta_k=\sup\{ \parallel \psi\parallel_2: \psi\in Z_k,\parallel \psi\parallel_{1,2}= 1\}$. We then have $\beta_k\rightarrow 0$ as $k\rightarrow \infty$.
\end{lemma}
Proof: By the definition of $\beta_k$, for each $j$ there exists $\psi_k\in Z_k$ such that $\parallel \psi_k\parallel_{1,2}= 1$  and
$\frac{1}{2}\beta_k<\parallel \psi_k\parallel_2$. According to the compactness of the embedding, we may assume (after taking a subsequence if necessary)that $\psi_k\rightharpoonup \psi$ weakly in $H^1(M,\Sigma M)$ and $\psi_k\rightarrow \psi$ strongly in $L^p(M,\Sigma M)$ for some
$\psi\in H^1(M,\Sigma M) $, where $1\leq p<2^*$. Therefore,
$$
1=\parallel \psi_k\parallel_{1,2}^2= \parallel D\psi_k\parallel_2^2+\parallel \psi_k\parallel_2^2=(1+\lambda_k^2)\parallel \psi_k\parallel_2^2
$$
where $\lambda_k$ is the eigenvalue of $\psi_k$, and $\lambda_k\rightarrow \infty$ as $k\rightarrow \infty$. Obviously, we have $\parallel\psi\parallel_2=0$ and
$\beta_k\rightarrow 0$ as $k\rightarrow \infty$.

\begin{lemma}\label{lem:4.3}
There exists  $\rho_k> r_k>0$ such that
\begin{flalign*}
&(A_1): a_k:=\inf_{\substack{\psi\in Z_k\\ \parallel \psi\parallel_{1,2}= r_k}}\mathfrak{L}(\psi)\rightarrow \infty, k\rightarrow \infty; &\\
&(A_2): b_k:=\sup_{\substack{\psi\in Y_k\\ \parallel \psi\parallel_{1,2}=\rho_k}}\mathfrak{L}(\psi)\leq 0\
\end{flalign*}
\end{lemma}
Proof: (i) Let $\psi\in Z_k$ with $\parallel \psi\parallel_{1,2}=r_k$, Then by $(\bf{H_1})$ and Lemma \ref{lem:4.2}, we have
\begin{eqnarray}\label{e:4.1}
 \mathfrak{L}(\psi)&=&\frac{1}{2}\parallel D\psi\parallel_2^2- \int_M H(x,\psi)dx\nonumber\\
  &\geq&\ \frac{1}{2}\parallel \psi\parallel_{1,2}^2-\frac{1}{2}\parallel \psi\parallel_2^2 -C\parallel \psi\parallel_{q}^q-C\nonumber\\
 &\geq&\ \frac{1}{2}r_k^2-\frac{1}{2}r_k^2\beta_k^2-Cr_k^q\beta_k^q-C
\end{eqnarray}
By Lemma \ref{lem:4.2}, we have $\beta_k^2\leq \frac{1}{2}$ for enough $k$, we obtain from (\ref{e:4.1}) that
\begin{equation}\label{e:4.2}
 \mathfrak{L}(\psi)\geq\ \frac{1}{4}r_k^2-Cr_k^q\beta_k^q-C
\end{equation}
As in the proof of [19,Lemma 3.6], choosing $r_k=(2pC\beta_k^q)^\frac{-1}{q-2}$, (\ref{e:4.2}) imply that
$$
\mathfrak{L}(\psi)\geq\ \frac{1}{4}r_k^2-Cr_k^q\beta_k^q-C\geq (\frac{1}{4}-\frac{1}{2p})(2qC\beta_k^q)^\frac{-2}{q-2}\rightarrow\infty
$$
as $k\rightarrow \infty$ and the condition $(A_1)$ is satisfied.
\par
(ii)Since $dim (Y_k)<\infty$, the $\parallel\cdot\parallel_2$ and $\parallel\cdot\parallel_{1,2}$-norms are equivalent on  $Y_k$, for any $\psi\in Y_k$,
there exist $C_k>\frac{1}{2}$ such that
\begin{equation}\label{e:4.3}
  \frac{1}{2}\int_M\mid D\psi\mid dx\leq \frac{1}{2}\parallel\psi\parallel_{1,2}^2\leq C_k\parallel\psi\parallel_2^2.
\end{equation}
By $(\bf{H_2})$, these exist $R_k>0$ such that
\begin{equation}\label{e:4.4}
  H(x,\psi)\geq 2C_k\mid \psi\mid^2,
\end{equation}
for any $(x,\psi)\in Y_k$ with $\mid \psi\mid\geq R_k$.

Choosing $M_k=\max\{\mid H(x,\psi)\mid:\mid \psi\mid\leq R_k, x\in M\}$, we obtain
\begin{equation}\label{e:4.5}
  H(x,\psi)\geq -M_k\geq 2C_k\mid\psi\mid^2-2C_kR_k^2-M_k
\end{equation}
for any $\mid \psi\mid\leq R_k$.

Using (\ref{e:4.4}) and (\ref{e:4.5}) we have
\begin{equation}\label{e:4.6}
  H(x,\psi)\geq 2C_k\mid\psi\mid^2-\overline{M_k}
\end{equation}
where $\overline{M_k}=2C_kR_k^2+M_k>0$.
Hence by (\ref{e:4.3}) and (\ref{e:4.6}), we obtain
$$
\mathfrak{L}(\psi)=\frac{1}{2}\int_M\mid D\psi\mid dx-\int_M H(x,\psi)dx\leq -C_k\parallel\psi\parallel_2^2+\overline{M_k}\mid M\mid.
$$
Therefore, we have $b_k\leq 0$ for enough $\rho_k>0$ where $\rho_k>r_k$. the condition $(A_2)$ is satisfied.

$\bf{Proof \ of \ Theorem1}$ According to Lemma \ref{lem:3.1}, $\mathfrak{L}$ satisfy the Cerami condition, Lemma \ref{lem:4.2} and Lemma \ref{lem:4.3} imply that the condition $(A_1)$ and $(A_2)$ is satisfied. By Fountain theorem there exists critical points $\psi_n\in H^1(M,\Sigma M)$ of $\mathfrak{L}$ such that $\mathfrak{L}(\psi_n)\rightarrow \infty$ as $n\rightarrow\infty$. Since $D^2\psi_n=D(D\psi_n)=H_\psi(x,\psi_n)$, where $H$ is $\alpha$-H\"{o}lder continuous in the direction of the base, where $0<\alpha<1$. By the Interior Schauder estimates in \cite{Amm},
we have $D\psi_n\in C^1(M,\Sigma M)$. Using the Interior Schauder estimates again, we obtain
those weak solutions $\psi_n$ are in fact $C^2(M,\Sigma M)$ and  classical solutions to problem (\ref{e:1.1}).
\par
$\bf{Founding}$ The *corresponding author: Xu Yang was supported by the NSFC
(grant no.11801499) of China.

%

\end{document}